\documentclass[12pt]{article}
\usepackage{amssymb}

\usepackage{amsmath,latexsym,amsfonts}
\usepackage{graphicx}
\usepackage{srcltx}
\usepackage{amscd}

\newtheorem{lemma}{Lemma}[section]
\newtheorem{coro}[lemma]{Corollary}
\newtheorem{prop}[lemma]{Proposition}
\newtheorem{thm}[lemma]{Theorem}
\newtheorem{defn}[lemma]{Definition}
\makeatletter\@addtoreset{equation}{section}
\renewcommand\theequation{\thesection.\@arabic\c@equation}

\begin{document}
\begin{center}
{\LARGE  Variation formulas for transversally harmonic and bi-harmonic maps }

 \renewcommand{\thefootnote}{}
\footnote{2000 \textit {Mathematics Subject Classification.}
53C12, 58E20}\footnote{\textit{Key words and phrases.}
Transversally harmonic map, transversally bi-harmonic map, trasnversal Jacobi operator,  variation formulas.
}
\renewcommand{\thefootnote}{\arabic{footnote}}
\setcounter{footnote}{0}

\vspace{0.5 cm} {\large Seoung Dal Jung }
\end{center}
\vspace{0.5cm}
\noindent {\bf Abstract.} In this paper, we study variation formulas for  transversally harmonic maps and bi-harmonic maps, respectively. We also study the transversal Jacobi field along a map and give several relations with infinitesimal automorphisms. 
\section{Introduction}
Let $(M,\mathcal F)$ and $(M',\mathcal F')$ be two foliated Riemannian manifolds and let $\phi:M\to M'$ be a smooth foliated map, i.e., $\phi$ is a leaf-preserving map. Then $\phi$ is
transversally harmonic if $\phi$ is a critical point of the transversal energy functional on any compact domain of $M$, which is defined in Section 3 (cf. [\ref{JJ2},\ref{KW1},\ref{KW2}]). Equivalently,  it is a solution of $\tau_b(\phi)=0$, where $\tau_b(\phi)$ is a transversal tension field, which is given by $\tau_b(\phi)={\rm tr}_Q\tilde\nabla d_T\phi$ (see [\ref{KW1}] for more details). That is, transversally harmonic maps are considered as harmonic maps between the leaf spaces [\ref{KW1},\ref{KW2}]. For harmonic maps, see [\ref{ES},\ref{XI}]. Also, we study the transversally bi-harmonic map as the critical point of the transversal bi-energy functional on any compact domain of $M$ (Section 6). In this paper, we study the second variation formulas for the transversal energy and transversal bi-energy of $\phi$. And we give some applications. This paper is organized as follows. In Section 2, we recall the basic facts on foliated manifolds. In Section 3, we review transversally harmonic maps and  the first variation formula. In Section 4, we give the second variation formula for the transversal energy. In Section 5, we define the transversal Jacobi operator along the foliated map  and study its realtion with infinitesimal automorphisms.  In Section 6, we study transversally bi-harmonic maps and their applications. In Section 7, we give  the second variation formula for the transversal bi-energy. Note that some results in Section 6 and Section 7 of the present paper can be found in [\ref{CW}], but the approach is different in a technical sense.
Throughout this paper,  $(M,\mathcal F)$ is considered as a foliated Riemannian manifold, i.e., a Riemannian manifold with a Riemannian foliation, and all leaves of $\mathcal F$ are compact.
\section{Preliminaries }
Let $(M,g,\mathcal F)$ be a $(p+q)$-dimensional  foliated Riemannian
manifold with foliation $\mathcal F$ of codimension $q$ and a bundle-like metric $g$ with respect to $\mathcal F$. Let $TM$ be the tangent bundle of $M$, $L$
the tangent bundle of $\mathcal F$, and  $Q=TM/L$ the
corresponding normal bundle of $\mathcal F$. 
  Let $g_Q$ be the holonomy invariant metric on
$Q$ induced by $g$. We denote by $\nabla^Q$ the transverse Levi-Civita
connection on the normal bundle $Q$ [\ref{TO1},\ref{TO2}]. Let $R^Q, K^Q,{\rm Ric}^Q $ and
$\sigma^Q$ be the transversal curvature tensor, transversal
sectional curvature, transversal Ricci operator and transversal
scalar curvature with respect to $\nabla^Q\equiv\nabla$, respectively.  Let $\Omega_B^r(\mathcal F)$ be the space of all {\it basic
$r$-forms}, i.e.,  $\omega\in\Omega_B^r(\mathcal F)$ if and only if
$i(X)\omega=0=i(X)d\omega$ for any $X\in\Gamma L$, where
$i(X)$ is the interior product. Then $\Omega^r(M)=\Omega_B^r(\mathcal F)\oplus \Omega_B^r(\mathcal F)^\perp$ [\ref{LO}]. Let $\kappa_B$ be the basic part of $\kappa$, the mean curvature form of $\mathcal F$. Then $\kappa_B$ is closed, i.e., $d\kappa_B=0$ [\ref{LO}]. The {\it basic
Laplacian} $\Delta_B$ acting on $\Omega_B^*(\mathcal F)$ is
defined by
\begin{equation}\label{2-1}
\Delta_B=d_B\delta_B+\delta_B d_B,
\end{equation}
where $\delta_B$ is the formal adjoint of
$d_B=d|_{\Omega_B^*(\mathcal F)}$[\ref{PR},\ref{TO2}].  Let $V(\mathcal F)$ be the space of all transversal infinitesimal automorphisms $Y$ of $\mathcal F$, i.e.,  $[Y,Z]\in \Gamma L$ for
all $Z\in \Gamma L$. Let $\bar V(\mathcal F)=\{\bar Y=\pi(Y)| Y\in V(\mathcal F)\}$, where $\pi:TM\to Q$ is a projection.
Trivially, $\bar V(\mathcal F)\cong \Omega_B^1(\mathcal F)$ [\ref{MO}].
For later use, we recall the transversal divergence theorem
[\ref{YO}] on a foliated Riemannian
manifold.
\begin{thm} \label{thm1-1}
Let $(M,g_M,\mathcal F)$ be a closed, oriented Riemannian manifold
with a transversally oriented foliation $\mathcal F$ and a
bundle-like metric $g_M$ with respect to $\mathcal F$. Then
\begin{equation}\label{2-4}
\int_M \operatorname{div_\nabla}\bar X = \int_M g_Q(\bar
X,\kappa_B^\sharp)
\end{equation}
for all $X\in V(\mathcal F)$, where $\operatorname{div_\nabla}\bar X$
denotes the transversal divergence of $\bar X$ with respect to the
connection $\nabla^Q$.
\end{thm}
Now, we define the bundle map $A_Y:\Lambda^r
Q^*\to\Lambda^r Q^*$ for any $Y\in V(\mathcal F)$ [\ref{KT}]
by
\begin{align}\label{2-5}
A_Y\omega =\theta(Y)\omega-\nabla_Y\omega \quad\forall\omega\in\Lambda^r Q^*,
\end{align}
where $\theta(Y)$ is the transverse Lie derivative. It is well-known [\ref{KT}] that, on $\Gamma Q$
\begin{align}\label{2-5-1}
A_Y s = -\nabla_{Y_s}\bar Y \quad\forall s\in\Gamma Q,
\end{align}
where $Y_s$ is the vector field such that $\pi(Y_s)=s$. So $A_Y$ depends only on $\bar Y=\pi(Y)$.
Since
$\theta(X)\omega=\nabla_X\omega$ for any $X\in\Gamma L$, $A_Y$
preserves the basic forms and depends only on $\bar Y$. 
   Let $E\to M$ be a vector bundle over $M$ and $\Omega_B^r (E)\equiv  \Omega_B^r(\mathcal F)\otimes E$ be the space of all $E$-valued basic $r$-forms. Let $\nabla$ be also the connection on $E$. Then the operator $A_X$ is extended to $\Omega_B^r(E)$ [\ref{JJ2}]. Now we define $d_\nabla :\Omega_B^r(E)\to \Omega_B^{r+1}(E)$ by
\begin{align}\label{2-6}
d_\nabla( \omega\otimes s)=  d_B \omega\otimes s+ (-1)^r\omega\wedge \nabla s
\end{align}
for any $\omega\in\Omega_B^r(\mathcal F)$ and $s\in E$.
Let $\delta_\nabla$ be the formal adjoint of $d_\nabla$. Then we define the Laplacian $\Delta$ on $\Omega_B^r(E)$  by
\begin{align}\label{2-6-1}
  \Delta=d_\nabla \delta_\nabla +\delta_\nabla d_\nabla.
  \end{align}
 From now on, let $\{E_a\}(a=1,\cdots,q)$ be a local orthonormal frame on $Q$ and $\theta^a$ be the $g_Q$-dual 1-form to $E_a$.
Then the  generalized Weitzenb\"ock type formula on $\Omega_B^r(E)$ is given by [\ref{JJ2}]
\begin{align}\label{2-7}
\Delta \Phi = \nabla_{\rm tr}^*\nabla_{\rm tr} \Phi
+ F(\Phi) + A_{\kappa_{B}^\sharp} \Phi,\quad\forall \Phi\in\Omega_B^r(E),
\end{align}
where  $\nabla_{\rm tr}^*\nabla_{\rm tr} =-\sum_a \nabla^2_{E_a,E_a}
+\nabla_{\kappa_B^\sharp}$ and  $F=\sum_{a,b=1}^{q}\theta^a\wedge i(E_b) R^\nabla(E_b,E_a)$. From (\ref{2-7}), we also have
\begin{align}\label{2-8}
\frac12\Delta_B| \Phi |^2 &= \langle\Delta \Phi, \Phi\rangle -
 |\nabla_{\rm tr} \Phi|^2 -
\langle A_{\kappa_{B}^\sharp}\Phi, \Phi\rangle  -\langle F(\Phi),\Phi\rangle,
\end{align}
where $\langle\cdot,\cdot\rangle$ is an inner product on $\Omega_B^r(E)$.

 Now, we recall the following generalized maximum principles.
\begin{lemma} $[\ref{JLK}]$ Let  $\mathcal F$ be a Riemannian foliation on a closed, oriented Riemannian manifold $(M,g_M)$. If $(\Delta_B -\kappa_B^\sharp)f\geq 0$ $($or $\leq 0)$ for any basic function $f$, then $f$ is constant.
\end{lemma}

\section{Transversally harmonic maps}
Let $(M,  g,\mathcal F)$  and $(M', g',\mathcal F')$  be
two foliated Riemannian manifolds and all leaves of $\mathcal F$ are compact.
  Let $\nabla^{M}$ and $\nabla^{M'}$ be the Levi-Civita connections on
$M$ and $M'$, respectively. And $\nabla$ and $\nabla'$ be the transverse
Levi-Civita connections on $Q$ and $Q'$, respectively.  Let $\phi:(M,g,\mathcal
F)\to (M', g',\mathcal F')$ be a smooth  foliated map,
i.e., $d\phi(L)\subset L'$.  We define $d_T\phi:Q \to Q'$  by
\begin{align}\label{3-1}
d_T\phi := \pi' \circ d \phi \circ \sigma,
\end{align}
where $\sigma : Q \to L^\perp$ is a bundle map satisfying $\pi\circ\sigma={\rm id}$. Then $d_T\phi$ is a section in $ Q^*\otimes
\phi^{-1}Q'$. Let $\nabla^\phi$
and $\tilde \nabla$ be the connections on $\phi^{-1}Q'$ and
$Q^*\otimes \phi^{-1}Q'$, respectively. Then $\phi:(M,g,\mathcal F)\to (M',g',\mathcal F')$ is called {\it transversally totally geodesic} if it satisfies
\begin{align}\label{3-2}
\tilde\nabla_{\rm tr}d_T\phi=0,
\end{align}
where $(\tilde\nabla_{\rm tr}d_T\phi)(X,Y)=(\tilde\nabla_X d_T\phi)(Y)$ for any $X,Y\in \Gamma Q$. And the {\it
transversal tension field} of $\phi$ is defined by
\begin{align}\label{3-3}
\tau_b(\phi)={\rm tr}_{Q}\tilde \nabla d_T
\phi=\sum_{a=1}^{q}(\tilde\nabla_{E_a} d_T\phi)(E_a),
\end{align}
where $\{E_a\}(a=1,\cdots,q)$ is a local orthonormal frame on $Q$.
Trivially, the transversal tension field $\tau_b(\phi)$ is a
section of $\phi^{-1}Q'$.
\begin{defn}{\rm
Let $\phi: (M, g,\mathcal F) \to (M', g',\mathcal F')$ be
a smooth foliated map. Then $\phi$ is said to be }
transversally harmonic {\rm if the transversal tension field
vanishes}, i.e., $\tau_b(\phi)=0$. 
\end{defn}
Let ${ vol}_L:M \to [0,\infty]$ be the volume map for which ${
vol}_L (x)$ is the volume of the leaf  passing through $x\in
M$. It is trivial that $vol_L$ is a basic function. And it holds [\ref{JJ2}] that
\begin{align}
d_B vol_L + (vol_L) \kappa_B =0.
\end{align}
 The transversal energy of $\phi$ on a compact domain $\Omega\subset
M$ is defined by
\begin{align}\label{eq2-4}
E_B(\phi;\Omega)={1\over 2}\int_{\Omega} | d_T \phi|^2 {1\over vol_L}
\mu_{M},
\end{align}\label{3-4}
where $\mu_{M}$ is the volume element of $M$.
Let
$V\in\phi^{-1}Q'$. Obviously, $V$ may be considered as a vector
field on $Q'$ along $\phi$. Then there is a 1-parameter family of
foliated maps $\phi_t$ with $\phi_0=\phi$ and ${d\phi_t\over
dt}|_{t=0}=V$. The family $\{\phi_t\}$ is said to be a {\it foliated variation} of $\phi$ with the normal variation vector field $V$. Then we have the first variation formula.

\begin{thm} $[\ref{JJ2}]$ \label{Thm3-2} {\rm (The first variation formula)} Let $\phi:(M,g,\mathcal F)\to (M',g',\mathcal F')$
be a smooth foliated map, and all leaves of $\mathcal F$ be compact. Let $\{\phi_t\}$ be a smooth foliated variation of $\phi$ supported in a compact domain $\Omega$. Then
\begin{align}\label{eq2-5}
{d\over dt}E_B(\phi_t,\Omega)|_{t=0}=-\int_{\Omega} \langle
V,\tau_b(\phi)\rangle {1\over vol_L}\mu_{M},
\end{align}
where $V(x)={d\phi_t\over dt}(x)|_{t=0}$ is the normal variation
vector field of $\{\phi_t\}$.
\end{thm}
\begin{defn}
 {\rm Let $\phi:(M,g,\mathcal F)\to (M',g',\mathcal F')$ be a smooth foliated map. Then the}  transversal stress-energy tensor $S_T(\phi)$ {\rm of $\phi$ is defined by}
 \begin{align}\label{3-6}
 S_T(\phi)={1\over 2}|d_T\phi|^2 g_{Q}-\phi^* g_{Q'},
 \end{align}
where $\phi^*$ is the pull-back of $\phi$.
 \end{defn}
 Trivially, $S_T(\phi)\in\otimes^2 Q^*$ is the symmetric 2-covariant normal tensor field on $M$.
 \begin{prop}$[\ref{CW}]$ Let $\phi:(M,g,\mathcal F)\to (M',g',\mathcal F')$ be a smooth foliated map. Then, for any vector field $X\in \Gamma Q$,
 \begin{align}
 ({\rm div}_\nabla S_T(\phi))(X)=-\langle\tau_b(\phi),d_T\phi(X)\rangle,
 \end{align}
 where $({\rm div}_\nabla S_T(\phi))(\cdot)=\sum_{a=1}^q (\nabla_{E_a}S_T(\phi))(E_a,\cdot)$.
 \end{prop}
 {\bf Proof.} Note that $[d_T\phi(X),d_T\phi(Y)]=d_T\phi([X,Y])$ for any  $X,Y\in\Gamma Q$. So, by direct calculation, the proof follows.  $\Box$

If ${\rm div}_\nabla S_T(\phi)=0$, then we say that $\phi$ satisfies the {\it transverse conservation law} [\ref{CW}]. 
The foliated map satisfying the transverse conservation law is said to be {\it transversally relatively harmonic}. Then we have the following.
\begin{coro} Any transversally harmonic map is  transversally relatively harmonic.
\end{coro}
 The converse of Corollary 3.5 does not hold. For the converse, see Theorem 7.4 below.

\noindent{\bf Remark.} $[\ref{JJ2}]$ Let $\phi:(M,g,\mathcal F)\to (M',g',\mathcal F')$ be a smooth foliated map. Then
\begin{align}\label{3-8}
d_\nabla d_T\phi=0,\quad \tilde\delta d_T\phi=-\tau_b(\phi),
\end{align}
where $\tilde\delta = \delta_\nabla -i(\kappa_B^\sharp)$.

\section{The second variation formula for the transversal energy}
Let $(M,  g,\mathcal F)$  and $(M', g',\mathcal F')$  be
two foliated Riemannian manifolds, and all  leaves of $\mathcal F$ be compact.
 Let $\phi:(M,g,\mathcal F)\to (M',g',\mathcal F')$ be a transversally harmonic map. For any $V,W\in \phi^{-1}Q'$, there exists a family of  foliated maps $\phi_{t,s} (-\epsilon<s,t<\epsilon)$ satisfying
 \begin{equation}\label{4-1}
\left\{\begin{split}&V={\partial\phi_{t,s}\over \partial t}\Big|_{(t,s)=(0,0)},\\
&W={\partial\phi_{t,s}\over \partial s}\Big|_{(t,s)=(0,0)},\\
&\phi_{0,0}=\phi.\qquad\qquad\qquad{}
\end{split}
\right.
\end{equation}
Then $\{\phi_{t,s}\}$ is said to be the {\it foliated variation} of $\phi$ with the {\it normal variation vector fields} $V$ and $W$.
Then we have the second variation formula for the transversal energy.
\begin{thm} {\rm (The second variation formula)}  Let $\phi:(M,g,\mathcal F)\to (M',g',\mathcal F')$ be a transversally harmonic map with $M$ compact without boundary, and all leaves of $\mathcal F$ be compact. Let $\{\phi_{t,s}\}$ be the foliated variation of $\phi$ with the normal variation vector fields $V$ and $W$. Then
\begin{align}\label{4-2}
&{\partial^2\over\partial t\partial s}E_B(\phi_{t,s})\Big|_{(t,s)=(0,0)}\\
&=\int_M \langle (\nabla_{\rm tr}^\phi)^*(\nabla_{\rm tr}^\phi) V-\nabla_{\kappa_B^\sharp}^\phi V-{\rm tr}_Q R^{Q'}(V,d_T\phi)d_T\phi,W\rangle {1\over vol_L}\mu_M,\notag
\end{align}
where ${\rm tr}_Q R^{Q'}(V,d_T\phi)d_T\phi=\sum_{a=1}^q R^{Q'}(V,d_T\phi(E_a))d_T\phi(E_a)$.
\end{thm}
{\bf Proof.}
  Let $\Phi:M \times (-\epsilon,\epsilon)\times(-\epsilon,\epsilon)\to M'$ be a smooth map, which is defined by $\Phi(x,t,s)=\phi_{t,s}(x)$. Let $\nabla^\Phi$ be the pull-back connection on $\Phi^{-1}Q'$. It is trivial that
$[X,{\partial\over\partial t}]=[X,{\partial\over\partial s}]=0$ for any vector field $X\in TM$. From the first normal variation formula (Theorem \ref{Thm3-2}), we have
\begin{align}\label{4-3}
{\partial\over\partial s}E_B(\phi_{t,s})=-\int_M\langle {\partial\Phi\over \partial s},\tau_b(\phi_{t,s})\rangle {1\over vol_L}\mu_M.
\end{align}
By differentiating (\ref{4-3}) with respect to $t$, we have
\begin{align*}
{\partial^2\over\partial t\partial s}E_B(\phi_{t,s})=-\int_M\{\langle {\partial^2\Phi\over\partial t\partial s},\tau_b(\phi_{t,s})\rangle +\langle{\partial\Phi\over\partial s},\nabla_{\partial\over\partial t}^\Phi\tau_b(\phi_{t,s})\rangle\}{1\over vol_L}\mu_M.
\end{align*}
At $(t,s)=(0,0)$, the first term vanishes since $\tau_b(\phi)=0$. Hence we  have
\begin{align}\label{4-4}
{\partial^2\over\partial t\partial s}E_B(\phi_{t,s})\Big |_{(t,s)=(0,0)}=-\int_M \langle W,\nabla^\Phi_{\partial\over\partial t}\tau_b(\phi_{t,s})\Big|_{(t,s)=(0,0)}\rangle{1\over vol_L}\mu_M.
\end{align}
We choose a local orthonormal basic frame field $\{E_a\}$ with $(\nabla E_a)(x)=0$. Then, at $x\in M$,
\begin{align*}
\nabla^\Phi_{\partial\over\partial t}\tau_b(\phi_{t,s})
&=\sum_a \{\nabla^\Phi_{\partial\over\partial t}\nabla^\Phi_{E_a}d_T\Phi(E_a)-\nabla^\Phi_{\partial\over\partial t} d_T\Phi(\nabla_{E_a}E_a)\}\\
&=\sum_a\{\nabla^\Phi_{E_a}\nabla^\Phi_{\partial\over\partial t} d_T\Phi(E_a)+R^\Phi({\partial\over\partial t},E_a)d_T\Phi(E_a)-\nabla^\Phi_{\partial\over\partial t} d_T\Phi(\nabla_{E_a}E_a)\}\\
&=\sum_a\{\nabla^\Phi_{E_a}\nabla^\Phi_{E_a}{\partial \Phi\over\partial t}-\nabla^\Phi_{\nabla_{E_a}E_a}{\partial \Phi\over\partial t}+R^{Q'}({\partial \Phi\over\partial t},d_T\Phi(E_a))d_T\Phi(E_a)\}.
\end{align*}
Hence, at $(t,s)=(0,0)$, we have
\begin{align}\label{4-5}
 &\nabla^\Phi_{\partial\over\partial t}\tau_b(\phi_{t,s})\Big|_{(t,s)=(0,0)}\notag\\
 &= \sum_a\{\nabla^\phi_{E_a}\nabla^\phi_{E_a}V-\nabla^\phi_{\nabla_{E_a}E_a}V +R^{Q'}(V,d_T\phi(E_a))d_T\phi(E_a)\}.
 \end{align}
 Hence the proof is complete. $\Box$
 
Let $\phi:(M,g,\mathcal F)\to (M',g',\mathcal F')$ be a foliated map with $M$ compact. Then we define the {\it transversal Hessian} $THess_\phi$ of $\phi$ by
\begin{align}\label{4-6}
THess_\phi(V,W)={\partial^2 \over\partial t \partial s}E_B(\phi_{t,s})\Big|_{(t,s)=(0,0)},
\end{align}
where $\{\phi_{t,s}\}$ is a foliated variation of $\phi$ with the normal variation vector fields $V$ and $W$. Then we have the following corollary.
 \begin{coro} Let $\phi:(M,g,\mathcal F)\to (M',g',\mathcal F')$ be a transversally harmonic map with $M$ compact without boundary, and all leaves of $\mathcal F$ be compact. Then, for any $V,W\in\phi^{-1}Q'$,
 \begin{align}\label{4-7}
 THess_\phi(V,W)&=\int_M\langle\nabla_{\rm tr}^\phi V,\nabla_{\rm tr}^\phi W\rangle{1\over vol_L}\mu_M\\
 &-\int_M\langle {\rm tr}_QR^{Q'}(V,d_T\phi)d_T\phi,W\rangle{1\over vol_L}\mu_M\notag
 \end{align}
 and $THess_\phi$ is symmetric, i.e., $THess_\phi(V,W)=THess_\phi(W,V)$ for any normal vector fields $V$ and $W$ along $\phi$.
 \end{coro}
 {\bf Proof.} From (\ref{4-2}) and (\ref{4-6}), we have
 \begin{align*}
& THess_\phi(V,W)\\
&=\sum_a\int_M\langle \nabla_{E_a}^\phi V,\nabla_{E_a}^\phi({1\over vol_L}W)\rangle\mu_M
-\int_M \langle\nabla_{(vol_L^{-1})\kappa_B^\sharp}^\phi V,W\rangle\mu_M \\
&-\int_M\langle{\rm tr}_QR^{Q'}(V,d_T\phi)d_T\phi,W\rangle{1\over vol_L}\mu_M\\
&=\int_M \sum_a\langle \nabla_{E_a}^\phi V,\nabla_{E_a}^\phi W\rangle{1\over vol_L}\mu_M -\int_M\langle{\rm tr}_QR^{Q'}(V,d_T\phi)d_T\phi,W\rangle{1\over vol_L}\mu_M\\
&-\int_M \langle \nabla^\phi_{d_B vol_L^\sharp + (vol_L)\kappa_B^\sharp}V,W\rangle {1\over vol_L^2}\mu_M.
 \end{align*}
  By (3.4), the last term in the last equality above vanishes. So the proof is completed. $\Box$

  If the transversal Hessian of $\phi:(M,g,\mathcal F)\to (M',g',\mathcal F')$ is positive semi-definite, i.e.,  $THess_\phi(V,V)\geq 0$  for any normal vector field $V$ along $\phi$, then $\phi$ is said to be {\it transversally stable}. From (\ref{4-7}), we have the following corollary.
  \begin{coro} {\rm (Stability)} Any transversally harmonic map from a compact(without boundary) foliated Riemannian manifold to  a foliated Riemannian manifold of non-positive transversal sectional curvature is transversally stable.
  \end{coro}
  \section{Transversal Jacobi operator along a map}
  Let $(M,g,\mathcal F)$ be a compact foliated Riemannian manifold, and all leaves of $\mathcal F$ are compact.
\begin{defn} {\rm Let $\phi:(M,g,\mathcal F)\to (M',g',\mathcal F')$ be a foliated map. Then the} transversal Jacobi operator $J_\phi^T:\Gamma \phi^{-1}Q'\to \Gamma \phi^{-1}Q'$ {\rm along $\phi$ is defined by
\begin{align}\label{5-1}
J^T_\phi(V)=(\nabla_{\rm tr}^\phi)^*(\nabla_{\rm tr}^\phi)V-\nabla_{\kappa_B^\sharp}^\phi V-{\rm tr}_Q R^{Q'}(V,d_T\phi)d_T\phi.
\end{align}
Any $V\in {\rm Ker}J^T_\phi$ is called a} transversal Jacobi field {\rm along $\phi$ for the transversal energy.}
\end{defn}
From (\ref{2-5-1}), the transversal Jacobi operator $J_{\nabla}^T\equiv J_{\rm id}^T$ along the identity map is given by
\begin{align}\label{5-1-1}
J_{\nabla}^T (\bar Y)=\nabla_{\rm tr}^*\nabla_{\rm tr}\bar Y -\rho^\nabla(\bar Y)+A_Y \kappa_B^\sharp
\end{align}
for any $Y\in V(\mathcal F)$, which is called to {\it generalized Jacobi operator} of $\mathcal F$ on $M$.  From (\ref{4-2}) and (\ref{4-6}), if $M$ is compact without boundary, then we have
\begin{align}\label{5-2}
THess_\phi(V,W)=\int_M\langle J^T_\phi(V),W\rangle {1\over vol_L}\mu_M.
\end{align}
Let $\{\phi_t\}$ be a smooth foliated variation of $\phi$ with the normal variation vector field $V$. From (4.4) and (5.3), we have
\begin{align}\label{5-3}
J^T_\phi(V)=-\nabla_{d\over dt} \tau_b(\phi_t)\Big|_{t=0}.
\end{align}
Hence we have the following proposition.
\begin{prop} Let $\phi:(M,g,\mathcal F)\to (M',g',\mathcal F')$ be a transversally harmonic map and $\{\phi_t\}$ be a smooth foliated variation of $\phi$ with the normal variation vector field $V$. Then $J^T_\phi(V)=0$, i.e., $V$ is a transversal Jacobi field along $\phi$ for the energy.
\end{prop}
Let $Y\in V(\mathcal F)$ be an infinitesimal automorphism on $(M,\mathcal F)$. It is well-known [\ref{JJ1}] that, if $\bar Y$ is a transversal affine field, i.e., $\theta(Y)\nabla=0$, then
\begin{align}\label{5-4}
\nabla_{\rm tr}^*\nabla_{\rm tr}\bar Y-\rho^\nabla(\bar Y)+A_Y\kappa_B^\sharp=0.
\end{align}
Hence we have the following proposition.
\begin{prop} On $(M,\mathcal F)$, any transversal affine  field is  a generalized Jacobi field for $\mathcal F$.
\end{prop}
{\bf Remark.} Since any transversal Killing field is transversal affine [\ref{KT}], any transversal Killing field is also generalized Jacobi field of $\mathcal F$. And the converse in Proposition 5.3 does not hold unless $\mathcal F$ is harmonic.
\begin{coro} $(cf. [\ref{JJ1}])$ On $(M,\mathcal F)$ with $M$ compact without boundary, the followings are equivalent: for any vector field $Y\in V(\mathcal F)$,

$(1)$ $\bar Y$ is a transversal Killing field, that is, $\theta(Y)g_Q=0$;

$(2)$ $\bar Y$ is a generalized Jacobi field for $\mathcal F$ satisfying $(i)$ ${\rm div}_\nabla(\bar Y)=0$ and $(ii)$ $\int_M g_Q(B_Y\bar Y,\kappa_B^\sharp)\geq 0$, where $B_Y=A_Y +A_Y^t$.
\end{coro}
Now, we recall the {\it transversal Jacobi operator} $J_\nabla:\Gamma Q\to \Gamma Q$ of $\mathcal F$ [\ref{KT}] by
\begin{align}
J_\nabla =\nabla_{\rm tr}^*\nabla_{\rm tr} -\rho^\nabla.
\end{align}
Then $\bar Y\in {\rm Ker}J_\nabla$ is  called a {\it transversal Jacobi field  for $\mathcal F$}      [\ref{KT}]. 
 Note that two operators $J_\nabla$ and $J_{\nabla}^T$ are related by
\begin{align}
  J_{\nabla}^T(\bar Y)=J_\nabla(\bar Y)+A_Y\kappa_B^\sharp.
  \end{align}  On a harmonic foliation $\mathcal F$, $J_{\nabla}^T =J_\nabla$. Hence, from  Proposition 5.3 and Corollary 5.4, we have the following corollary.
\begin{coro} $[\ref{KTT}]$ Let $\mathcal F$ be a harmonic foliation on a compact Riemannian manifold $(M,g)$. Then the following are equivalent:

 $(1)$ $\bar Y$ is a transversal Killing field, i.e., $\theta(Y)g_Q=0$;
 
  $(2)$ $\bar Y$ is a transversal Jacobi field of $\mathcal F$ and ${\rm div}_\nabla(\bar Y)=0$;
  
   $(3)$ $\bar Y$ is transversal affine field, i.e., $\theta(Y)\nabla=0$. 
\end{coro}
Now, we have the vanishing theorem about the transversal Jacobi field along the map.
\begin{thm} Let $(M,g,\mathcal F)$ be a closed,connected Riemannian manifold with a foliation $\mathcal F$ and a bundle-like metric $g$. Assume the transversal Ricci operator is non-positive and negative at some point. Then any generalized Jacobi field $\bar Y$ for $\mathcal F$ is trivial, i.e, $Y$ is tangential to $\mathcal F$.
\end{thm}
{\bf Proof.} It is well-known [\ref{JJ1}] that
\begin{align}\label{5-6}
\frac12(\Delta_B-\kappa_B^\sharp)|\bar Y|^2 = g_Q(J_{\nabla}^T(\bar Y)+g_Q(\rho^\nabla(\bar Y),\bar Y)-|\nabla_{\rm tr}\bar Y|^2.
\end{align} 
Let $\bar Y$ be a generalized Jacobi field for $\mathcal F$. Then we have
\begin{align}\label{5-7}
\frac12(\Delta_B-\kappa_B^\sharp)|\bar Y|^2 = g_Q(\rho^\nabla(\bar Y),\bar Y)-|\nabla_{\rm tr}\bar Y|^2.
\end{align}
Since the transversal Ricci curvature is non-positive, we have
\begin{align}\label{5-8}
(\Delta_B-\kappa_B^\sharp)|\bar Y|^2\leq 0.
\end{align}
Hence, by the generalized maximum principle (Lemma 2.2), $|\bar Y|$ is constant. Again, from (\ref{5-7}), $\bar Y$ is parallel. Moreover, since $\rho^\nabla$ is negative at some point, $\bar Y$ is trivial. Equivalently, $Y$ is tangential to $\mathcal F$. $\Box$

If $\bar Y$ is a transversal Jacobi field of $\mathcal F$, i.e., $J_\nabla(\bar Y)=0$, then $J_{\nabla}^T(\bar Y)=-\nabla_{\kappa_B^\sharp}\bar Y$.  Therefore,   from (\ref{5-6}) we have
  \begin{align}\label{5-12}
 \frac12 \Delta_B |\bar Y|^2 =g_Q(\rho^\nabla(\bar Y),\bar Y) -|\nabla_{\rm tr}\bar Y|^2.
  \end{align}
 From (\ref{5-12}),  we have the following corollary [\ref{KT}, p535].
\begin{coro} $[\ref{KT}]$ Let $(M,g,\mathcal F)$ be as in Theorem 5.6. Assume the transversal Ricci operator is non-positive and negative at some point. Then any transversal Jacobi field of $\mathcal F$ is trivial, i.e., $Y$ is tangential to $\mathcal F$.
\end{coro}
{\bf Proof.} Let $\bar Y$ be a transversal Jacobi field of $\mathcal F$. Since $\rho^\nabla$ is non-positive, from (\ref{5-12}), $\Delta_B |\bar Y|^2\leq 0$. Since $\Delta=\Delta_B$ on a basic function, by the maximum principle, $|\bar Y|$ is constant. Since $\rho^\nabla$ is negative at some point, from (5.11), $\bar Y$ is trivial. $\Box$

  \section{Transversally bi-harmonic maps }

Let $(M,  g,\mathcal F)$  and $(M', g',\mathcal F')$  be
two foliated Riemannian manifolds.
   Let $\phi:(M,g,\mathcal
F)\to (M', g',\mathcal F')$ be a smooth  foliated map. 
 Now we define the {\it transversal bi-tension field $(\tau_2)_b(\phi)$ } of $\phi$ by
\begin{align}\label{6-9}
(\tau_2)_b(\phi)=J_\phi^T(\tau_b(\phi)).
\end{align}
\begin{defn} {\rm Let $\phi:(M,g,\mathcal F)\to (M',g',\mathcal F')$ be a smooth foliated map. Then $\phi$ is said to be} transversally bi-harmonic {\rm if the transversal bi-tension field vanishes, i.e., $(\tau_2)_b(\phi)=0$.}
\end{defn}
Trivially, $\phi$ is a transversally bi-harmonic map if and only if the transversal tension field $\tau_b(\phi)$ is a transversal Jacobi field along $\phi$. Moreover, any transversal harmonic map is a transversal bi-harmonic map.
Now, we define the {\it transversal bi-energy} of $\phi$ supported in a compact domain $\Omega$ by
\begin{align}\label{6-12}
(E_2)_B(\phi,\Omega)=\frac12\int_\Omega |\tilde\delta d_T\phi|^2{1\over vol_L}\mu_M.
\end{align}
Then we have the following theorem.
\begin{thm} $(${\rm The first variation formula for the transversal bi-energy}$)$ Let $\phi:(M,g,\mathcal F)\to (M',g',\mathcal F')$ be a smooth foliated map, and all leaves of $\mathcal F$ be compact. Let $\{\phi_t\}$ be a foliated variation of $\phi$ with the variation vector field $V$ in a compact domain $\Omega$.
Then we have
\begin{align}\label{6-13}
{d\over dt}(E_2)_B(\phi_t,\Omega)\Big|_{t=0}=-\int_\Omega \langle (\tau_2)_b(\phi),V\rangle{1\over vol_L}\mu_M.
\end{align}
\end{thm}
{\bf Proof.} Let $\{\phi_t\}$ be a foliated variation of $\phi$ such that ${d\phi_t\over dt}\Big|_{t=0}=V$ and $\phi_0=\phi$. Choose a local orthonormal basic frame $\{E_a\}$ with $(\nabla E_a)(x)=0$. Define $\Phi:M\times (-\epsilon,\epsilon)\to M'$ by $\Phi(x,t)=\phi_t(x)$. Let $\nabla^\Phi$ be the pull-back connection on $\Phi^{-1}Q'$. Obviously, $d_T\Phi(E_a)=d_T\phi(E_a)$ and $d\Phi({d\over dt})={d\phi_t\over dt}$. Moreover, it is trivial that $\nabla_{d\over dt}{d\over dt}=\nabla_{d\over dt}E_a =\nabla_{E_a}{d\over dt}=0$. Hence, from (\ref{3-8}), we have
\begin{align}\label{4-4-1}
{d\over dt}(E_2)_B(\phi_t,\Omega)&=\int_\Omega \langle \nabla_{d\over dt}^\Phi \tau_b(\phi_t),\tau_b(\phi_t)\rangle{1\over vol_L}\mu_M.
\end{align}
From (\ref{5-3}), it follows that
\begin{align*}
{d\over dt}(E_2)_B(\phi_t,\Omega)|_{t=0}&=-\int_\Omega \langle J_\phi^T(V),\tau_b(\phi)\rangle{1\over vol_L}\mu_M\\
&=-\int_\Omega\langle J_\phi^T(\tau_b(\phi)),V\rangle{1\over vol_L}\mu_M.
\end{align*}
The last equality above follows from (5.3) and the symmetry of the transversal Hessian $THess$ of $\phi$.
From (\ref{6-9}), the proof is complete. $\Box$
\begin{coro} Let $\phi:(M,g,\mathcal F)\to (M',g',\mathcal F')$ be a smooth foliated map, and all leaves of $\mathcal F$ be compact. Then $\phi$ is transversally bi-harmonic if and only if it is a critical point of the transversal bi-energy $(E_2)_B(\phi)$ of $\phi$ on any compact domain.
\end{coro}
Then we have the following (cf. [\ref{CW}]).
\begin{thm} Let $\phi:(M,g,\mathcal F)\to (M',g',\mathcal F')$ be a transversally bi-harmonic map with $M$ compact without boundary, and all leaves of $\mathcal F$ be compact. Assume that the transversal sectional curvature $K^{Q'}$ of $\mathcal F'$ is non-positive. Then $\phi$ is transversally harmonic.
\end{thm}
{\bf Proof.} Let $\{E_a\}$ be a local orthonormal basic frame of $Q$. Then $(\tau_2)_b(\phi)=0$ implies that
\begin{align}\label{6-14}
(\nabla_{\rm tr}^\phi)^*\nabla_{\rm tr}^\phi \tau_b(\phi)-\nabla_{\kappa_B^\sharp}^\phi \tau_b(\phi)-\sum_a R^{Q'}(\tau_b(\phi),d_T\phi(E_a))d_T\phi(E_a)=0.
\end{align}
From (\ref{2-7}) and (\ref{2-8}), we have $\Delta_B |\tau_b(\phi)|^2 =2\langle (\nabla_{\rm tr}^\phi)^*\nabla_{\rm tr}^\phi \tau_b(\phi),\tau_b(\phi)\rangle -2|\nabla_{\rm tr}\tau_b(\phi)|^2$. Hence from (\ref{6-14}), we have
\begin{align}\label{6-15}
{1\over 2}(\Delta_B-\kappa_B^\sharp)|\tau_b(\phi)|^2 &=-|\nabla_{\rm tr}\tau_b(\phi)|^2 \notag\\&
+ \sum_a\langle R^{Q'}(\tau_b(\phi),d_T\phi(E_a))d_T\phi(E_a),\tau_b(\phi)\rangle.
  \end{align}
 Since the transversal sectional curvature $K^{Q'}$ of $\mathcal F'$ is non-positive, we have
 \begin{align*}
 (\Delta_B-\kappa_B^\sharp)|\tau_b(\phi)|^2\leq 0.
 \end{align*}
 Hence, by the generalized maximum principle (Lemma 2.3), $|\tau_b(\phi)|$ is constant.
Again, from (\ref{6-15}), we have that for all $a$,
\begin{align}\label{eq5-7}
\nabla_{E_a}\tau_b(\phi)=0.
\end{align}
Now, we define the normal vector bundle $X$ by
\begin{align*}
X={1\over vol_L}\sum_a\langle d_T\phi(E_a),\tau_b(\phi)\rangle E_a.
\end{align*}
Then we have
\begin{align*}
{\rm div}_\nabla (X)&=-{1\over {vol_L}^{2}}\langle d_T\phi(d_B vol_L^\sharp),\tau_b(\phi)\rangle +{1\over vol_L}
|\tau_b(\phi)|^2\\
&={1\over vol_L}\langle d_T\phi(\kappa_B^\sharp),\tau_b(\phi)\rangle +{1\over vol_L}|\tau_b(\phi)|^2.
\end{align*}
The last equality above follows from Lemma 2.2. By integrating and by using the transversal divergence theorem (Theorem 2.1), we have
\begin{align}\label{6-16}
\int_M |\tau_b(\phi)|^2{1\over vol_L}\mu_M=0,
\end{align}
which implies that $\tau_b(\phi)=0$. So $\phi$ is transversally harmonic.  $\Box$

\section{The second variation formula for the transversal bi-energy}
Let $(M,g,\mathcal F)$ and $(M',g',\mathcal F')$ be two foliated Riemannian manifolds. Let $\phi:(M,g,\mathcal F)\to (M',g',\mathcal F')$ be a transversally bi-harmonic map. Then we have the second variational formula of the transversal bi-energy as follows.
\begin{thm}  $($ {\rm The second variation formula for the transversal bi-energy}$)$ Let $\phi:(M,g,\mathcal F)\to (M',g',\mathcal F')$ be a smooth foliated map with $M$ compact without boundary and all leaves be compact. Let $\{\phi_t\}$ be a foliated variation of $\phi$ with the normal variation vector field $V$. Then
\begin{align*}
&{d^2\over dt^2} (E_2)_B (\phi_t)\Big|_{t=0}\\
&=\int_M \{|(\tau_2)_b(\phi))|^2 -\langle (\tau_2)_b(\phi)),\nabla_VV\rangle-\langle R^{Q'}(V,\tau_b(\phi))\tau_b(\phi),V\rangle\}{1\over vol_L}\mu_M\\
&-2\sum_a\int_M\langle(\nabla_{E_a}R^{Q'})(V,d_T\phi(E_a))\tau_b(\phi),V\rangle{1\over vol_L}\mu_M\\
&+\sum_a\int_M\langle (\nabla_{\tau_b(\phi)}R^{Q'})(V,d_T\phi(E_a))d_T\phi(E_a),V\rangle{1\over vol_L}\mu_M\\
&-4\sum_a\int_M\langle R^{Q'}(\nabla^\phi_{E_a}V,\tau_b(\phi))d_T\phi(E_a),V\rangle{1\over vol_L}\mu_M.
\end{align*}
\end{thm}
{\bf Proof.} Let $V\in \phi^{-1}Q'$ and let $\{\phi_t\}$ be a foliated variation of $\phi$ such that ${d\phi_t\over dt}\Big|_{t=0}=V$ and $\phi_0=\phi$. We choose a local orthonormal basic frame $\{E_a\}$ such that $(\nabla E_a)(x)=0$ at point $x$. Define $\Phi:M\times (-\epsilon,\epsilon)\to M'$ by $\Phi(x,t)=\phi_t(x)$. Let $\nabla^\Phi$ be the pull-back connection on $\Phi^{-1}Q'$. Obviously, $d_T\Phi(E_a)=d_T\phi(E_a)$ and ${\partial\Phi\over \partial t}={d\phi_t\over dt}$. Trivially, $\nabla_{d\over dt}{d\over dt}=\nabla_{d\over dt}E_a =\nabla_{E_a}{d\over dt}=0$. Hence, from (\ref{4-4-1}) we have
\begin{align}\label{7-1}
{d^2\over dt^2}(E_2)_B(\phi_t)&=\int_M \Big(\langle\nabla^\Phi_{d\over dt} \nabla^\Phi_{d\over dt}\tau_b(\phi_t),\tau_b(\phi_t)\rangle+|\nabla^\Phi_{d\over dt} \tau_b(\phi_t)|^2\Big){1\over vol_L}\mu_M.
\end{align}
From (\ref{4-5}), we have
\begin{align*}
\nabla^\Phi_{d\over dt} \tau_b(\phi_t)=\sum_a(\nabla^\Phi)^2_{E_a,E_a}{\partial\Phi\over \partial t} + \sum_aR^{Q'}({\partial\Phi\over \partial t},d_T\phi_t(E_a))d_T\phi_t(E_a),
\end{align*}
and then
\begin{align*}
\nabla^\Phi_{d\over dt}\nabla^\Phi_{d\over dt} \tau_b(\phi_t)=\sum_a\nabla^\Phi_{d\over dt}\Big((\nabla^\Phi)^2_{E_a,E_a} {\partial\Phi\over \partial t} + R^{Q'}({\partial\Phi\over \partial t},d_T\phi_t(E_a))d_T\phi_t(E_a)\Big).
\end{align*}
By a long calculation, we get
\begin{align*}
\nabla_{d\over dt}^\Phi(\nabla^\Phi)^2_{E_a,E_a}{\partial\Phi\over \partial t}&=(\nabla^\Phi)^2_{E_a,E_a}\nabla^\Phi_{d\over dt}{\partial\Phi\over \partial t}+\nabla^\Phi_{E_a}R^\Phi({d\over dt},E_a){\partial\Phi\over \partial t}\\
&+R^\Phi({d\over dt},E_a)\nabla^\Phi_{E_a}{\partial\Phi\over \partial t} +R^\Phi(\nabla_{E_a}E_a,{d\over dt}){\partial\Phi\over \partial t}
\end{align*}
and
\begin{align*}
&\nabla^\Phi_{d\over dt} R^{Q'}({\partial\Phi\over \partial t},d_T\phi_t(E_a))d_T\phi_t(E_a)\\
&=(\nabla_{\partial\Phi\over \partial t}R^{Q'})({\partial\Phi\over \partial t},d_T\phi_t(E_a))d_T\phi_t(E_a)+R^{Q'}(\nabla^\Phi_{d\over dt}{\partial\Phi\over \partial t},d_T\phi_t(E_a))d_T\phi_t(E_a)\\
&+R^{Q'}({\partial\Phi\over \partial t},\nabla^\Phi_{d\over dt}d_T\phi_t(E_a))d_T\phi_t(E_a)+R^{Q'}({\partial\Phi\over\partial t},d_T\phi_t(E_a))\nabla^\Phi_{d\over dt}d_T\phi_t(E_a).
\end{align*}
Since $[{d\over dt},E_a]=0$, we have $\nabla_V d_T\phi(E_a)=\nabla_{d_T\phi(E_a)}V=\nabla^\phi_{E_a}V$. Hence from the equations above, we have
\begin{align*}
&\nabla^\Phi_{d\over dt}\nabla^\Phi_{d\over dt} \tau_b(\phi_t)|_{t=0}\\
&=-J^T_\phi(\nabla_VV)+\sum_a\nabla^\phi_{E_a}R^{Q'}(V,d_T\phi(E_a))V\\
&+2\sum_aR^{Q'}(V,d_T\phi(E_a))\nabla^\phi_{E_a}V +\sum_aR^{Q'}(d_T\phi(\nabla_{E_a}E_a),V)V\\
&+\sum_a(\nabla_V R^{Q'})(V,d_T\phi(E_a))d_T\phi(E_a)+\sum_aR^{Q'}(V,\nabla^\phi_{E_a}V)d_T\phi(E_a).
\end{align*}
 So, by the first and second Bianchi identites, we have
\begin{align*}
\langle\nabla^\Phi_{d\over dt} \nabla_{d\over dt}^\Phi \tau_b(\phi_t),\tau_b(\phi_t)\rangle|_{t=0}
&=-\langle J_\phi^T(\nabla_VV),\tau_b(\phi)\rangle +\langle R^{Q'}(V,\tau_b(\phi))V,\tau_b(\phi)\rangle \\
&+\sum_a\langle (\nabla_V R^{Q'})(V,d_T\phi(E_a))d_T\phi(E_a),\tau_b(\phi)\rangle\\
&+\sum_a\langle (\nabla_{E_a}R^{Q'})(V,d_T\phi(E_a))V,\tau_b(\phi)\rangle\\
& +4\sum_a\langle R^{Q'}(V,d_T\phi(E_a))\nabla^\phi_{E_a}V,\tau_b(\phi)\rangle.
\end{align*}
By integrating together with (\ref{5-3}), we have
\begin{align*}
&\int_M\langle\nabla^\Phi_{d\over dt}\nabla^\Phi_{d\over dt} \tau_b(\phi_t)|_{t=0},\tau_b(\phi)\rangle{1\over vol_L}\mu_M\\
&=-\int_M\langle J_\phi^T(\tau_b(\phi)),\nabla_VV\rangle{1\over vol_L}\mu_M-\int_M\langle R^{Q'}(V,\tau_b(\phi))\tau_b(\phi),V\rangle{1\over vol_L}\mu_M\notag\\
&-2\sum_a\int_M\langle(\nabla_{E_a}R^{Q'})(V,d_T\phi(E_a))\tau_b(\phi),V\rangle{1\over vol_L}\mu_M \notag\\
&+\sum_a\int_M\langle(\nabla_{\tau_b(\phi)} R^{Q'})(V,d_T\phi(E_a))d_T\phi(E_a),V\rangle{1\over vol_L}\mu_M\notag\\
& -4\sum_a\int_M\langle R^{Q'}(\nabla^\phi_{E_a}V,\tau_b(\phi))d_T\phi(E_a),V\rangle{1\over vol_L}\mu_M.
\end{align*}
From (\ref{6-9}) and (\ref{7-1}), the proof follows. $\Box$

Then we have the following corollary (cf. [\ref{CW}]).
\begin{coro} Let $\phi:(M,g,\mathcal F)\to (M',g',\mathcal F')$ be a transversally bi-harmonic map with $M$ compact without boundary. Let $\{\phi_t\}$ be a foliated variation of $\phi$ with the normal variation vector field $V$. Then
\begin{align*}
{d^2\over dt^2} (E_2)_B (\phi_t)\Big|_{t=0}=&-\int_M\langle R^{Q'}(V,\tau_b(\phi))\tau_b(\phi),V\rangle{1\over vol_L}\mu_M\\
&-2\sum_a\int_M\langle(\nabla_{E_a}R^{Q'})(V,d_T\phi(E_a))\tau_b(\phi),V\rangle{1\over vol_L}\mu_M\\
&+\sum_a\int_M\langle (\nabla_{\tau_b(\phi)}R^{Q'})(V,d_T\phi(E_a))d_T\phi(E_a),V\rangle{1\over vol_L}\mu_M\\
&-4\sum_a\int_M\langle R^{Q'}(\nabla_{E_a}^\phi V,\tau_b(\phi))d_T\phi(E_a),V\rangle{1\over vol_L}\mu_M.
\end{align*}
\end{coro}
\begin{defn} {\rm If the transversally bi-harmonic map $\phi:(M,g,\mathcal F)\to (M',g',\mathcal F')$ satisfies ${d^2\over dt^2}(E_2)_B(\phi_t)\Big|_{t=0}\geq 0$, then $\phi$ is said to be} stable.
\end{defn}
Note that any transversally harmonic map can be considered as a stable transversally bi-harmonic map. We also have the following theorem ([\ref{CW}]).
\begin{thm} Let $(M,\mathcal F)$ be a closed Riemannian manifold with a foliation $\mathcal F$ and let $(M',\mathcal F')$ be a Riemannian manifold with a constant transversal sectional curvature $C>0$. If a foliated map $\phi:(M,\mathcal F)\to (M',\mathcal F')$ is  stable transversally bi-harmonic and relatively harmonic, then $\phi$ is transversally harmonic.
\end{thm}
{\bf Proof.} Since the transversal sectionsal curvature $K^{Q'}$ of $\mathcal F'$ is constant $C>0$, from Corollary 7.2, we have
\begin{align*}
{d^2\over dt^2} (E_2)_B (\phi_t)\Big|_{t=0}=&-\int_M\langle R^{Q'}(V,\tau_b(\phi))\tau_b(\phi),V\rangle{1\over vol_L}\mu_M\\
&-4\sum_a\int_M\langle R^{Q'}(\nabla_{E_a}^\phi V,\tau_b(\phi))d_T\phi(E_a),V\rangle{1\over vol_L}\mu_M.
\end{align*}
Let $V=\tau_b(\phi)$. Since $\phi$ is relatively harmonic, i.e., $\langle \tau_b(\phi),d_T\phi(X)\rangle$ for any vector field $X\in \Gamma Q'$, we have
\begin{align*}
{d^2\over dt^2} (E_2)_B (\phi_t)\Big|_{t=0}&=-4\sum_a\int_M\langle R^{Q'}(\nabla_{E_a}^\phi \tau_b(\phi),\tau_b(\phi))d_T\phi(E_a),\tau_b(\phi)\rangle{1\over vol_L}\mu_M\\
&=4C\sum_a\int_M \langle \nabla_{E_a}^\phi \tau_b(\phi),d_B\phi(E_a)\rangle |\tau_b(\phi)|^2{1\over vol_L}\mu_M\\
&=-4C\int_M|\tau_b(\phi)|^4{1\over vol_L}\mu_M.
\end{align*}
This stability implies that $\tau_b(\phi)=0$, i.e., $\phi$ is transversally harmonic. $\Box$

\noindent{\bf Acknowledgements.} This research was supported by the Basic Science Research Program through the National Research Foundation of Korea(NRF) funded by the Ministry of Education, Science and Technology (2010-0021005) and NRF-2011-616-C00040.

\noindent Seoung Dal Jung

\noindent Department of Mathematics and Research Institute for Basic Sciences, Jeju National University,
Jeju 690-756, Korea

\noindent {\it E-mail address} : sdjung@jejunu.ac.kr




\end{document}